\documentclass[11pt]{amsart}

\usepackage{amsfonts,amssymb,amsthm,amscd,enumerate}

\newtheorem{theorem}{Theorem}

\title[Symmetric polynomials]{Symmetric polynomials and non-finitely generated $Sym (\mathbb N)$-invariant ideals}

\author{Eudes Antonio da Costa}

\address{Eudes Antonio da Costa: Departamento de Matem\'atica, Universidade de Bras\'\i lia, 70910-900 Bras\'\i lia, DF, Brasil}

\email{eudes@uft.edu.br}

\author{Alexei Krasilnikov}

\address{Alexei Krasilnikov: Departamento de Matem\'atica, Universidade de Bras\'\i lia, 70910-900 Bras\'\i lia, DF, Brasil}

\email{alexei@unb.br}

\date{}

\begin{document}

\begin{abstract}
Let $K$ be a field and let $\mathbb N = \{ 1,2, \dots \}$. Let $R_n=K[x_{ij} \mid 1\le i\le n, j\in \mathbb N]$ be the ring of polynomials in $x_{ij}$ $(1 \le i \le n, j \in \mathbb N )$ over $K$. Let $S_n = Sym (\{ 1,2, \ldots, n \})$ and $Sym (\mathbb N)$ be the groups of the permutations of the sets $\{ 1,2,\dots, n \}$ and $\mathbb N$, respectively. Then $S_n$ and $Sym (\mathbb N)$ act on $R_n$ in a natural way: $\tau (x_{ij})=x_{\tau(i)j}$ and $\sigma (x_{ij})=x_{i\sigma (j)}$ for all $\tau \in S_n$ and $\sigma \in Sym(\mathbb N)$. Let $\overline{R}_n$ be the subalgebra of the symmetric polynomials in $R_n$,
\[
\overline{R}_n = \{ f \in R_n \mid \tau (f) = f \mbox{ for each
} \tau  \in S_n \} .
\]
In 1992 the second author proved that if $char (K)= 0$ or $char(K)=p > n$ then every $Sym (\mathbb N)$-invariant ideal in $\overline{R}_n$ is finitely generated (as such). In this note we prove that this is not the case if $char (K)=p\le n$.

We also survey some results about $Sym (\mathbb N)$-invariant ideals in polynomial algebras and some related results.
\end{abstract}

\maketitle

\section{Introduction}

Let $K$ be a unital associative and commutative ring and let $\mathbb N = \{ 1,2, \dots \}$ be the set of all positive integers. Let $R_n=K[x_{ij} \mid 1\le i\le n, j\in \mathbb N]$ be the ring of polynomials in $x_{ij}$ $(1 \le i \le n, j \in \mathbb N )$ over $K$.

For a non-empty set $A$, let $Sym (A)$ denote the group of all permutations of $A$. The group $Sym (\mathbb N)$ acts on $R_n$ in a natural way: $\sigma (x_{ij}) = x_{i \sigma (j)}$ if $\sigma \in Sym (\mathbb N)$.  An ideal $I$ of $R_n$ is called \textit{ $Sym (\mathbb N)$-invariant} if $\sigma (I) = I$ for all $\sigma \in Sym( \mathbb N)$.

It is clear that $R_n$ contains ideals that are not finitely generated. However, the following theorem holds.

\begin{theorem}[see \cite{AH07, Cohen67, Cohen87, HS12}]
\label{theorem1}
Let $K$ be a Noetherian unital associative and commutative ring. Then each $Sym (\mathbb N)$-invariant ideal of $R_n=K[x_{ij} \mid 1\le i\le n, j\in \mathbb N]$ is finitely generated (as a $Sym (\mathbb N)$-invariant ideal).
\end{theorem}

For $n=1$ this theorem was proved by Cohen \cite{Cohen67} in 1967 and rediscovered independently by Aschenbrenner and Hillar \cite{AH07} in 2007; for an arbitrary positive $n$ this was proved by Cohen \cite{Cohen87} in 1987 and rediscovered independently by Hillar and Sullivant \cite{HS12} in 2012. Cohen's results were motivated by the finite basis problem for identities of metabelian groups and the results of Aschenbrenner, Hillar and Sullivant by applications to chemistry and algebraic statistics.

Let
\[
d_{i_1 i_2 \dots i_n} = det \left(
\begin{array}{cccc}
x_{1 i_1} & x_{1 i_2} & \dots & x_{1 i_n}
\\
\dots & \dots & \dots & \dots
\\
x_{n i_1} & x_{n i_2} & \dots & x_{ni_n}
\end{array}
\right)
\]
be the determinant of the matrix above. It is clear that $d_{i_1 i_2 \dots i_n}$ is a polynomial contained in $R_n$. Let $L_n$ be the subalgebra in $R_n$ generated by all polynomials $d_{i_1 i_2 \dots i_n}$ $(i_{\ell} \in \mathbb N)$. The group $Sym (\mathbb N)$ acts on $L_n$ in a natural way: $\sigma (d_{i_1 i_2 \dots i_n }) = d_{\sigma (i_1) \sigma (i_2) \dots \sigma (i_n)}$ ($\sigma \in Sym (\mathbb N)$). $Sym (\mathbb N)$-invariant ideals of $L_n$ are also defined in a natural way.

The following theorem was proved by Draisma \cite{Draisma10} in 2010; it solves a problem arising from applications to algebraic  statistics and chemistry posed in \cite{AH07}.

\begin{theorem}[see \cite{Draisma10}]
\label{D10}
Let $K$ be a field of characteristic $0$. Then every $Sym (\mathbb N)$-invariant ideal in $L_n$ is finitely generated (as such).
\end{theorem}

However, Theorem \ref{D10} is not valid over a field $K$ of characteristic $2$: the algebra $L_2$ over such $K$ contains $Sym (\mathbb N)$-invariant ideals that are not finitely generated. More precisely, the following theorem holds.

\begin{theorem}[see \cite{DK13}]
\label{DK13}
Suppose that $K$ is a field of characteristic $2$. Let $I$ be the $Sym( \mathbb N)$-invariant ideal in $L_2$ generated by the set
\[
\{ d_{12}d_{23} \dots d_{(k-1)k} d_{1k} \mid k = 3,4, \dots \} .
\]
Then $I$ is not finitely generated (as a $Sym (\mathbb N)$-invariant ideal in $L_2$).
\end{theorem}

The proof of Theorem \ref{DK13} is based on the ideas of Vaughan-Lee \cite{VL75} developed in order to construct an example of a non-finitely based variety of abelian-by-nilpotent Lie algebras.

Theorems \ref{D10} and \ref{DK13} raise the following problem:

\medskip
\noindent
\textbf{Problem.} \textit{For which $n$ and $p$ the algebra $L_n$ over a field $K$ of characteristic $p>0$ is $Sym (\mathbb N)$-Noetherian (that is, each $Sym (\mathbb N)$-invariant ideal in $L_n$ is finitely generated)?}

\medskip
In order to find an approach to this problem one can consider a simpler (but similar in a certain sense) question about $Sym (\mathbb N)$-Noetherianity of the subalgebra $\overline{R}_n$ of symmetric polynomials (defined below).

Let $S_n = Sym( \{ 1,2, \dots, n \} )$. The group $S_n$ acts on $R_n$ in a natural way: $\tau (x_{ij}) = x_{\tau (i) j}$. In other words, $S_n$ acts on the infinite matrix
\begin{equation}
\label{matrix}
\left(
\begin{array}{ccccc}
x_{11} & x_{12} & \dots & x_{1i} & \dots
\\
\dots & \dots & \dots & \dots & \dots
\\
x_{n1} & x_{n2} & \dots & x_{ni} & \dots
\end{array}
\right)
\end{equation}
permuting its lines.

We call a polynomial $f \in R_n$ \textit{symmetric} if $\tau (f) = f$ for all $\tau \in S_n$. Let $\overline{R}_n$ be the set of all symmetric polynomials of $R_n$; it is clear that $\overline{R}_n$ is a $Sym (\mathbb N)$-invariant $K$-subalgebra in $R_n$. $Sym (\mathbb N)$-invariant ideals of $\overline{R}_n$ are defined in a natural way.

The following theorem was proved by the second author of the present note \cite[Lemma 7]{Kras92} in 1992 in order to prove the finite basis property for certain varieties of nilpotent-by-abelian Lie algebras.

\begin{theorem}[see \cite{Kras92}]
\label{theorem_Kras92}
Let $K$ be a field of characteristic $0$ or of characteristic $p > n$. Then every $Sym (\mathbb N)$-invariant ideal in $\overline{R}_n$ is finitely generated (as such).
\end{theorem}

The proof of Theorem \ref{theorem_Kras92} is very simple. Let $\pi : R_n \rightarrow \overline{R}_n$ be the symmetrization map:
\begin{multline*}
\pi (f) = \pi \bigl( f(x_{11}, x_{12}, \dots ; \dots ; x_{n1}, x_{n2}, \dots ) \bigr)
\\
= \frac{1}{n!} \sum_{\tau \in S_n} f(x_{\tau (1) 1}, x_{\tau (1) 2}, \dots ; \dots ; x_{\tau (n) 1}, x_{\tau (n) 2}, \dots )
\end{multline*}
for all $f \in R_n$. Note that $\pi$ is $K$-linear, commutes with each permutation $\sigma \in Sym (\mathbb N)$ and if $\overline{a} \in \overline{R}_n$, $f \in R_n$ then $\pi (\overline{a}) = \overline{a}$ and $\pi (\overline{a} f) = \overline{a} \, \pi (f)$.

Let $I$ be a $Sym (\mathbb N)$-invariant ideal in $\overline{R}_n$. Then $I \cdot R_n$ is a $Sym (\mathbb N)$-invariant ideal in $R_n$ that, by Theorem \ref{theorem1}, is generated (as a $Sym (\mathbb N)$-invariant ideal in $R_n$) by a finite set, say, by $ b_1,  \dots , b_k $. We may assume that $b_i \in I \subseteq \overline{R}_n$ for all $i$.

We claim that $b_1, \dots , b_k$ generate $I$ as a $Sym (\mathbb N)$-invariant ideal in $\overline{R}_n$. Indeed, take an arbitrary element $g \in I$. Since $b_1,  \dots , b_k$ generate $I \cdot R_n$ as a $Sym (\mathbb N)$-invariant ideal in $R_n$, we have $g = \sum_{i} \sigma_i (b_{\ell_i}) f_i$ for some $\sigma_i \in Sym (\mathbb N)$ and $f_i \in R_n$. It follows that
\[
g = \pi (g) = \pi \bigl( \sum_{i} \sigma_i (b_{\ell_i}) f_i \bigr) = \sum_{i} \pi \bigl( \sigma_i (b_{\ell_i}) f_i \bigr) = \sum_{i} \sigma_i (b_{\ell_i}) \pi (f_i),
\]
that is, $g$ belongs to the $Sym (\mathbb N)$-invariant ideal of $\overline{R}_n$ generated by $b_1,  \dots , b_k$, as claimed.

The proof of Theorem \ref{D10} is much more sophisticated than that of Theorem \ref{theorem_Kras92} but follows the same plan, with the group $S_n$ acting on the matrix (\ref{matrix}) by line permutations replaced by the group $SL_n(K)$ acting on (\ref{matrix}) by multiplications from the left and with the symmetrization map replaced by the Reynolds operator.

The aim of the present note is to prove the following theorem.

\begin{theorem}
\label{R_n}
Let $K$ be a field of characteristic $p$ such that $0 < p \le n$. Then the algebra $\overline{R}_n$ contains $Sym (\mathbb N)$-invariant ideals that are not finitely generated.
\end{theorem}

Thus, the algebra $\overline{R}_n$ is $Sym (\mathbb N)$-Noetherian if $char K =0$ or $char K >n$ and is not $Sym (\mathbb N)$-Noetherian if $0 < char K \le n$. Possibly (although it still remains unknown), this is also the case for the algebra $L_n$ if $char K >2$.

Theorem \ref{R_n} is a corollary of the following result. For each $k \in \mathbb N$, let $h_k$ be the polynomial in $R_n$ defined as follows:
\[
h_k = x_{11} x_{12} \dots x_{1k} + x_{21} x_{22} \dots x_{2k} + \dots + x_{n1} x_{n2} \dots x_{nk} .
\]
For example, $h_1 =  x_{11} + x_{21}+ \dots + x_{n1}$, $h_2 = x_{11} x_{12} + x_{21}x_{22} + \dots +x_{n1}x_{n2}$. It is clear that $h_k \in \overline{R}_n$ for all $k$.

Let $U$ be the $Sym(\mathbb N)$-invariant ideal of $\overline{R}_n$ generated by the set
\[
\{ h_k \in \overline{R}_n \mid k = 1,2, \ldots \} .
\]
Our main result is as follows.

\begin{theorem}
\label{maintheorem_1}
Let $K$ be a field of characteristic $p>0$. Suppose that $n \ge p$. Then the ideal  $U$ is not finitely generated as a $Sym(\mathbb N)$-invariant ideal in $\overline{R}_n$.
\end{theorem}

We will prove Theorem \ref{maintheorem_1} by proving the following slightly stronger result:

\begin{theorem}
\label{maintheorem_2}
Let $K$ be a field of characteristic $p>0$. Suppose that $n \ge p$. Then, for each $k \in \mathbb N$, the polynomial $h_k$ is not contained in the $Sym (\mathbb N)$-invariant ideal of $\overline{R}_n$ generated by the set $\{ h_l \mid l \in \mathbb N, l \ne k \}$.
\end{theorem}

\noindent
\textbf{Remarks.} 1. For other recent results about $Sym (\mathbb N)$-Noetherian polynomial algebras see, for instance, articles \cite{BD11,DEKL13,DKutt09,DKutt13,HdC13} and a survey \cite{Draisma13}.

2. The main result of \cite{Kras92} about polynomial rings is, in fact, stronger then Theorem \ref{theorem_Kras92}. It has been proved there that if $char K =0$ or $char K >n$ then each $Sym (\mathbb N)$-invariant $\overline{R}_n$-submodule in $R_n$ is finitely generated (as a $Sym (\mathbb N)$-invariant $\overline{R}_n$-submodule).

3. Some results about Noetherian properties of polynomial rings in infinitely many variables have been proved in \cite{Kras89,Kras90} in order to solve the finite basis problem for certain varieties of groups and group representations.

Let $\Psi$ be the set of endomorphisms $\psi_{k \ell}$ $(k \ne \ell )$ of $R_n=K[x_{ij} \mid 1\le i\le n, j\in \mathbb N]$ such that
\[
\psi_{k \ell} (x_{ij}) =
\left\{
\begin{array}{cc}
x_{i k} x_{i \ell} & \mbox {if } j = \ell ;
\\
x_{ij} & \mbox{ otherwise}.
\end{array}
\right .
\]
We say that a subset $S$ of $R_n$ is $\Psi$-closed if $\psi_{k \ell} (S) \subseteq S$ for all $\psi_{k \ell} \in \Psi$. It has been proved in \cite{Kras89} that, over a Noetherian associative and commutative unital ring $K$, each $Sym (\mathbb N)$-invariant $\Psi$-closed $\overline{R}_n$-submodule in $R_n$ is finitely generated (as such). Moreover, over such $K$ each $Sym (\mathbb N)$-invariant $\Psi$-closed $\widetilde{R}_n$-submodule in $R_n$ is finitely
generated \cite{Kras90}. Here $\widetilde{R}_n$ is a $K$-subalgebra  of $\overline{R}_n$ generated by all products of the form
\[
f(x_{11}, x_{12}, \dots ) f(x_{21}, x_{22}, \dots ) \dots f(x_{n1}, x_{n2}, \dots )
\]
where $f(t_1, t_2, \dots ) \in K[t_i \mid i \in \mathbb N ]$.

Note that the result about $\overline{R}_n$-modules has been proved using techniques similar to one used by Cohen \cite{Cohen67} while the proof of the result about $\widetilde{R}_n$-modules is based on different ideas.

4. Let $M$ be the free metabelian group on a free generating set $\{ x_i \mid i \in \mathbb N \}$. Then the group $Sym (\mathbb N)$ acts on $M$ permuting the free generators $x_i$. Cohen \cite{Cohen67} has proved that in $M$ all $Sym (\mathbb N)$-invariant normal subgroups are finitely generated (as such). For an application of this result see, for example, \cite{DK03}.

Note that $Sym (\mathbb N)$-invariant normal subgroups are finitely generated (as $Sym (\mathbb N)$-inva\-riant normal subgroups) not just in the variety of metabelian groups but also in some larger varieties (Vaughan-Lee \cite{VL70}).  However, it can be easily shown that in the variety of centre-by-metabelian groups (and, therefore, in all larger varieties) this property does not hold.

\section{Proof of Theorem \ref{maintheorem_2} }

We claim that to prove Theorem \ref{maintheorem_2} one may assume without loss of generality that $n = p$. Indeed, suppose that $n>p$. Let $\psi : R_n \rightarrow R_p$ be the homomorphism of $R_n$ onto $R_p$ such that
\[
\psi (x_{ij}) =
\left\{
\begin{array}{ll}
 x_{ij} & \mbox{ if } 1 \le i \le p;
\\
0 & \mbox{ if } p < i \le n.
\end{array}
\right.
\]
It is clear that $\psi (\overline{R}_n) = \overline{R}_p$ and $\psi \sigma = \sigma \psi$ for all $\sigma \in Sym (\mathbb N)$. Hence, to prove that the polynomial
\[
h_k = x_{11} x_{12} \dots x_{1k} + x_{21} x_{22} \dots x_{2k} + \dots + x_{n1} x_{n2} \dots x_{nk}
\]
is not contained in the $Sym (\mathbb N)$-invariant ideal of $\overline{R}_n$ generated by the set $\{ h_{\ell} \mid l \in \mathbb N, \ell \ne k \}$ it suffices to prove that the polynomial
\[
\psi ( h_k ) = x_{11} x_{12} \dots x_{1k} + x_{21} x_{22} \dots x_{2k} + \dots + x_{p1} x_{p2} \dots x_{pk}
\]
is not contained in the $Sym (\mathbb N)$-invariant ideal of $\psi (\overline{R}_n) = \overline{R}_p$ generated by the set  $\{ \psi ( h_{\ell}) \mid l \in \mathbb N, \ell \ne k \}$ where
\[
\psi ( h_{\ell} ) = x_{11} x_{12} \dots x_{1 \ell} + x_{21} x_{22} \dots x_{2 \ell} + \dots + x_{p1} x_{p2} \dots x_{p \ell} .
\]
The claim follows.

Further, let $\mathbb F_p = \mathbb Z / p \mathbb Z$ be the prime subfield of $K$, $\mathbb F_p < K$. We claim that to prove Theorem \ref{maintheorem_2} one can assume without loss of generality that $K = \mathbb F_p$. Indeed, let $R_{p, {\mathbb F_p}} = \mathbb F_p [x_{ij} \mid 1 \le i \le p; j \in \mathbb N]$ and let $\overline{R}_{p, {\mathbb F_p}}$ be the subalgebra of symmetric polynomials of  $R_{p, {\mathbb F_p}}$. Then $R_{p, {\mathbb F_p}} < R_p = K [x_{ij} \mid 1 \le i \le p; j \in \mathbb N]$, $\overline{R}_{p, {\mathbb F_p}} < \overline{R}_p$ and $h_k \in \overline{R}_{p, {\mathbb F_p}}$ for all $k$.

Suppose that $h_k$ is contained in the $Sym (\mathbb N)$-invariant ideal of $\overline{R}_p$ generated by the set $\{ h_{\ell} \mid l \in \mathbb N, \ell \ne k \}$. Then
\[
h_k = \sum_{s} \sigma_s  (h_{\ell_s}) f_s
\]
where $\sigma_s \in Sym (\mathbb N)$, $f_s \in \overline{R}_p$, $\ell_s <k$ for all $s$. Let $\mathcal B$ be a basis of $K$ viewed as a vector space over $\mathbb F_p$ such that $1 \in \mathcal B$. Then, for each $s$, $f_s = f_{s,0} + \sum_t b_t f_{s,t}$ where $f_{s,0} \in R_{p, \mathbb F_p}$ and, for all $t$, $b_t \in \mathcal B \setminus \{ 1 \}$, $f_{s,t} \in R_{p, \mathbb F_p}$. It follows that
\[
h_k = \sum_{s} \sigma_s  (h_{\ell_s}) f_{s,0} + \sum_s \sum_t b_t \, \sigma_s  (h_{\ell_s}) f_{s,t}.
\]
Since $h_k$,  $\sigma_s  (h_{\ell_s}) f_{s,0}$, $ \sigma_s  (h_{\ell_s}) f_{s,t} \in R_{p, \mathbb F_p}$, we have
\begin{equation}
\label{hkhl}
h_k = \sum_{s} \sigma_s  (h_{\ell_s}) f_{s,0}.
\end{equation}
Note that $f_{s,0} \in \overline{R}_{p, \mathbb F_p}$ for all $s$. Indeed, for each $\tau \in S_n$,
\[
f_{s,0} + \sum_t b_t f_{s,t} = f_s = \tau (f_s) = \tau (f_{s,0}) + \sum_t b_t \tau (f_{s,t})
\]
so $\tau (f_{s,0} ) = f_{s,0}$. Thus, if $h_k$ is contained in the $Sym (\mathbb N)$-invariant ideal of $\overline{R}_p$ generated by the set $\{ h_{\ell} \mid l \in \mathbb N, \ell \ne k \}$ then, by (\ref{hkhl}), $h_k$ belongs to the $Sym (\mathbb N)$-invariant ideal of $\overline{R}_{p, \mathbb F_p}$ generated by $\{ h_{\ell} \mid l \in \mathbb N, \ell \ne k \}$. Hence, one may assume that $K = \mathbb F_p$, as claimed.

From now on, we assume $K = \mathbb F_p$ and write $R_p$ and $\overline{R}_p$ for $R_{p, \mathbb F_p}$ and $\overline{R}_{p, \mathbb F_p}$, respectively.

Let $R = \mathbb Z[x_{ij} \mid 1\le i\le p, j \in \mathbb N]$. Let $\mu : R \rightarrow R_{p}$ be the homomorphism of $R$ onto $R_p$ such that $\mu ( x_{ij}) = x_{ij}$ for all $i,j$. Suppose, in order to get a contradiction, that $h_k$ is contained in the $Sym (\mathbb N)$-invariant ideal of $\overline{R}_p$ generated by the set $\{ h_{\ell} \mid l \in \mathbb N, \ell \ne k \}$. Then in $\overline{R}_p$ we have
\[
h_k = \sum_{s} \sigma_s  (h_{\ell_s}) f_s
\]
where $\sigma_s \in Sym (\mathbb N)$, $f_s \in \overline{R}_p$, $\ell_s <k$ for all $s$. Hence, in $R$ we have
\begin{equation}
\label{hk}
h_k = \sum_{s} \sigma_s  (h_{\ell_s}) g_s + p \ g
\end{equation}
where $g_s, g \in R$, $\mu (g_s ) = f_s$ for all $s$. It is clear that we may assume that $g_s \in \overline{R}$ for all $s$. It follows that $g = \frac{1}{p} (h_k - \sum_{s} \sigma_s  (h_{\ell_s}) g_s) \in \overline{R}$.

Let
\[
m = x_{11}^{u_{11}} x_{21}^{u_{21}} \dots x_{p1}^{u_{p1}} \ x_{12}^{u_{12}} x_{22}^{u_{22}} \dots x_{p2}^{u_{p2}} \dots x_{1k}^{u_{1k}} x_{2k}^{u_{2k}} \dots x_{pk}^{u_{pk}} \in R
\]
be a monomial. Define a \textit{multi-degree} $d(m)$ of the monomial $m$ as follows:
\[
d(m)=(u_{1},u_2, \dots, u_{k}, 0,0, \dots),
\]
where $u_k = u_{1k} + u_{2k} + \dots + u_{pk}$ for all $k \in \mathbb N$. Note that $d(m \cdot m')=d(m)+d(m')$ for all monomials $m, m'$ in $R$. Note also that we may assume that the polynomials $h_k$, $\sigma_s  (h_{\ell_s}) g_s$ and $g$ in (\ref{hk}) are of the same multi-degree $(1,1, \dots , 1, 0,0, \dots )$.

Let $S = \mathbb Z [t_{j} \mid j \in \mathbb N ]$. Define a homomorphism $\eta : R \rightarrow S$ of $R$ onto $S$ by $\eta (x_{ij}) = t_j$. It is clear that $\eta (h_k) = p \ t_1 t_2 \dots t_k$, $\eta \bigl( \sigma_s (h_{\ell_s}) \bigr) = p \ t_{\sigma_s (1)} t_{\sigma_s (2)} \dots t_{\sigma_s (\ell_s )}$.

Note that if $m \in R$ is a monomial of multi-degree $(\varepsilon_1 ,\varepsilon_2, \dots , \varepsilon_k, 0, 0, \dots)$ where $\varepsilon_j = 0,1$ then the number of elements of the $S_p$-orbit of $m$ is a multiple of $p$. Hence, $\eta (m) = p \, (q m')$ for some monomial $m' \in S$ and some integer $q \in \mathbb Z$. It follows that, for all $s$, $\eta (g_s) = p \, u_s$  and $\eta (g) = p \, u$ for some $u_s, u \in S$.

Applying $\eta$ to the both sides of the equality (\ref{hk}), we get
\[
p \ t_1 t_2 \dots t_k = p^2 \ \sum_s t_{\sigma_s (1)} t_{\sigma_s (2)} \dots t_{\sigma_s (\ell_s )} u_s + p^2 \ u
\]
for some $u_s, u \in S$. This is a contradiction because the right hand side of the equality above is a (non-zero) multiple of $p^2$ in $S = \mathbb Z [t_j \mid j \in \mathbb N ]$ and the left hand side is not. It follows that $h_k$ is not contained in the $Sym (\mathbb N)$-invariant ideal of $\overline{R}_p$ generated by the set $\{ h_{\ell} \mid l \in \mathbb N, \ell \ne k \}$, as required.

The proof of Theorem \ref{maintheorem_2} is completed.

\bigskip \noindent
\textbf{Acknowledgements.}
The first author was supported by CNPq grant 554712/2009--1. The second author was partially supported by CNPq grant 307328/2012--0, by DPP/UnB and by RFBR grant 11--01--00945.

\end{document}